\documentclass[12pt, oneside]{amsart}

\usepackage{amsmath,amssymb}
\usepackage{accents}
\usepackage{graphicx}
\usepackage{verbatim}
\usepackage{color}
\usepackage{float}
\usepackage{enumerate}
\usepackage[textwidth=6in, textheight=9in,marginpar=0.75in]{geometry}

\newlength{\dhatheight}

\newtheorem{theorem}{Theorem}[section]

\newcommand{\nidt}{\noindent}

\newcommand{\del}{\partial}

\begin{document}
\title[The smooth 4-genus of the prime knots through 12 crossings]
{The smooth 4-genus of (the rest of) the prime knots through 12 crossings}

\author[M.~Brittenham]{Mark Brittenham}
\address{Department of Mathematics\\
        University of Nebraska\\
         Lincoln NE 68588-0130, USA}
\email{mbrittenham2@unl.edu}

\author[S.~Hermiller]{Susan Hermiller}
\address{Department of Mathematics\\
        University of Nebraska\\
         Lincoln NE 68588-0130, USA}
\email{hermiller@unl.edu}

\date{November 30, 2021}

\begin{abstract}
We compute the smooth 4-genera of the prime knots with 12 crossings
whose values, as reported on the KnotInfo website, were unknown.
This completes the calculation of the smooth 4-genus for all prime
knots with 12 or fewer crossings.
\end{abstract}

\maketitle


\section{Introduction}\label{sec:intro}


The {\it genus} (or {\it 3-genus}) $g(K)$ of a knot $K$ in the 3-sphere $S^3$ is the minimum genus 
of an orientable spanning surface (or {\it Seifert surface}) $F$ with $\del F=K$ embedded 
in $S^3$. This is one of the most basic measures of the complexity of a knot, and following 
the development of knot Floer homology (\cite{oz04},\cite{ras03}) and the efforts of many 
researchers (\cite{juh10},\cite{km97},\cite{oz04b}), determination of the 3-genus
for a given knot can now be viewed 
as a `routine' computation, in particular since the ability to compute knot Floer homology was recently 
added to the computer program SnapPy~\cite{snappy}.

On the other hand, the {\it smooth 4-genus} $g_4(K)$ (which we will refer
to as the {\it 4-genus} throughout this paper), which is the minimum genus of a 
smooth, orientable, properly embedded surface $\Sigma$ in the standard smooth 4-ball
$B^4$ with $K=\Sigma\cap\partial B^4$, currently enjoys no such advantage. There
are several computable lower bounds, both classical and modern, 
including the inequality
\begin{equation}\label{eq:lowerbound}
g_4(K) \ge |\sigma(K)|/2,
\end{equation}
where $\sigma(K)$ is the knot 
signature~\cite{mur65}. 
Rasmussen's $s$-invariant \cite{ras10} and Ozsv\'ath-Szab\'o's $\nu$-invariant \cite{oz11}
also provide lower bounds on the 4-genus. 
The classical invariants 3-genus $g(K)$ and unknotting number $u(K)$
give upper bounds for the 4-genus, but $g_4(K)=g(K)$ is relatively rare
(only 277 of the 2977 prime knots through 12 crossings have this equality),
and the inequality $g_4(K)\leq u(K)$ is more often used to help determine unknotting
number from a known value of $g_4(K)$, rather than the other way around. 
In practice, useful upper bounds on $g_4(K)$, in the end, typically come from direct construction
of a smooth orientable spanning surface $\Sigma$.

Despite these seeming limitations, among the 2977 prime knots with 12 or fewer crossings, the smooth
4-genus is known for all but 3 of them 
(listed in Theorem~\ref{theorem:slice}(a) below) at the time of this writing, 
according the KnotInfo website~\cite{knotinfo}. 
In this paper in Theorem~\ref{theorem:slice}(a)
we determine the value of $g_4(K)$ for these 3 remaining knots, completing
the computation of the 4-genus for all prime knots with crossing number at most 12.

For 2205 of the 2977 prime knots with crossing number at most 12, 
upper bounds have been computed to show that 
(combined with Equation~(\ref{eq:lowerbound})) the 4-genus satisfies the equality
 $g_4(K) = |\sigma(K)|/2$~\cite{knotinfo}.
In 2015 Lewark and McCoy~\cite{lewmcc19}
utilized genus one concordances (see Section~\ref{sec:genus1} for definitions
and background) to compute upper bounds for smooth 4-genera, and
used algebraic methods (discussed more in Section~\ref{sec:up} below)
to find lower bounds on 4-genera.
They determined the (previously unknown) smooth 4-genus of 639 knots and, 
as a result, 
there remained at the time only 22
knots of up to 12 crossings whose smooth 4-genera were still unknown. In the years 
since then, the 4-genus of the Conway knot $K11n34$ was determined in the spectacular result of Piccirillo~\cite{pic20} that the Conway knot is not 
smoothly slice, and more recently
(during the course of our work on this problem) Karageorghis and Swenton 
announced~\cite{karswe21} the computation
of the values of $g_4(K)$ for 18 of the remaining 21 knots 
(listed in Theorem~\ref{theorem:slice}(b) below). Their approach involved the 
calculation of the double slice genus $g_{ds}(K)$ of these knots (the minimum genus of an 
unknotted surface $\Sigma$ in $S^4$ with $\Sigma\cap S^3=K$, 
where $S^3$ is the equatorial 3-sphere); for these 18 knots they
find that
$g_{ds}(K)=2$ or $3$, implying that $g_4(K)\leq 1$.

Our approach is motivated by the work of Lewark and McCoy in~\cite{lewmcc19}.
For the three remaining knots for which the 4-genus was unknown,
we use a combination of algebraic methods together with extended methods for
a computer search for 
genus one concordances in order to determine the necessary lower
bounds on the 4-genus. 
We also expand on the techniques
of~\cite{lewmcc19} for finding upper bounds on 4-genera through genus one concordances
in order to directly construct a genus one spanning surface for 
the 18 knots for which the 4-genus was determined by
Karageorghis and Swenton, giving an alternate proof of that result. 
In all of our work, we restrict our view 
to four specific operations for constructing genus one
concordances: crossing change,
switching a pair of oppositely signed crossings, or resolving
or de-resolving a pair of crossings (see Section~\ref{sec:genus1} for 
more details).

\begin{theorem} \label{theorem:slice}
(a) Each of the 3 knots $K12a153$, $K12n239$ and $K12n512$ is genus 
one concordant, via de-resolution of a pair of crossings,
to a knot with 4-genus 3, and hence each has smooth 4-genus equal to~2. 

(b) Each of the 18 knots 
$K11n80$, $K12a187$, $K12a230$, $K12a317$, $K12a450$, $K12a570$, $K12a624$,
$K12a636$, $K12a905$, $K12a1189$, $K12a1208$, $K12n52$, $K12n63$, $K12n225$, $K12n555$, 
$K12n558$, $K12n665$ and $K12n886$ is genus one concordant, via crossing resolution
or switching two crossings, to one of
the slice knots $6_1$, $8_8$, $8_{20}$, $10_{75}$, $10_{87}$, $10_{137}$, $K11n74$, 
or $K12n256$, and hence each of the 18 knots has smooth 4-genus equal to~1.
\end{theorem}

We note that our computation of $g_4(K)=2$ for the
three knots $K$ in Theorem~\ref{theorem:slice}(a),
together with Karageorghis and Swenton's work~\cite{karswe21}
which shows that $g_{ds}(K)\leq 4$, also
determines that the double slice genus satisfies
$g_{ds}(K)=4$ for the knots $K$ in the list $K12a153$, $K12n239$ and $K12n512$.


\section{Genus one concordances}\label{sec:genus1}


A key tool in this paper is what should properly 
be thought of as the 4-genus version of a crossing change. A {\it genus one
concordance} between two knots $K_0$, $K_1$ is a twice-punctured torus 
$\Sigma\cong T^2\setminus\textrm{int}(D^2_1\coprod D^2_2)$ 
properly and smoothly embedded in $S^3\times[0,1]$ such that 
$\Sigma\cap (S^3\times\{0\})=K_0$ and $\Sigma\cap (S^3\times\{1\})=K_1$.
Gluing $\Sigma\subseteq S^3\times[0,1]$ to any $g_4$-minimizing surface in $B^4$
bounding $K_1$ yields a smooth surface in $B^4$ of genus
$g_4(K_1)+1$ with boundary $K_0$, and so $g_4(K_0)\leq g_4(K_1)+1$. Switching
roles of $K_0$ and $K_1$ and using the same genus one concordance shows that 
$g_4(K_1)\leq g_4(K_0)+1$, and so 
\begin{equation}\label{eq:concordance}
g_4(K_1)-1\leq g_4(K_0)\leq g_4(K_1)+1
\end{equation}
when $K_0$ and $K_1$ are genus one concordant.

Our goal is to use these inequalities to 
determine $g_4(K_0)$. We will do this in two ways.
The first is to find a genus one concordance to a knot $K_1$
for which it is known that $g_4(K_0)\geq g_4(K_1)+1$
(in our applications, $g_4(K_1)=0$), 
and so (with Equation~(\ref{eq:concordance})) 
$g_4(K_0)=g_4(K_1)+1$. In this case we say that the genus one concordance {\it pulls
the 4-genus of $K_0$ down} to $g_4(K_1)+1$.
The second is to find a genus one concordance to a knot $K_1$
for which it is known that $g_4(K_0)\leq g_4(K_1)-1$
(in our applications,  $g_4(K_1)=3$), 
and so $g_4(K_0)=g_4(K_1)-1$. In this case we say that 
the genus one concordance {\it pulls
the 4-genus of $K_0$ up} to $g_4(K_1)-1$.
In particlar, since all of the 21 knots
that we wish to consider here have 4-genus either 1 or 2 (from KnotInfo~\cite{knotinfo}), we seek a
genus one concordance to a knot $K_1$ with $g_4(K_1)=0$ (i.e., a smoothly slice knot),
`pulling down' 
to establish $g_4(K_0)=1$, or to a knot $K_1$ with $g_4(K_1)=3$, `pulling up' to 
establish $g_4(K_0)=2$. 

\medskip

\begin{figure}[h]
\begin{center}
\includegraphics[width=3.5in]{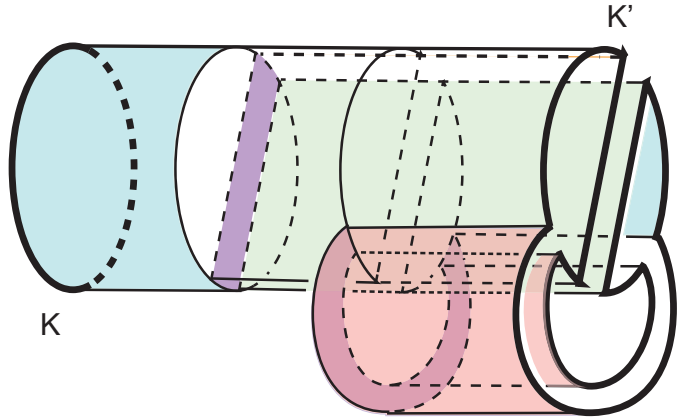}
\caption{A genus one concordance: schematic}\label{fig:genus1}
\end{center}
\end{figure}

We employ four standard operations for building genus one
concordances from a knot $K$
in the process of making most of our 4-genus determinations.
Each of these operations can
can be carried out on the level of knot diagrams, and  
they are especially amenable to implementation in
software when $K$ is given as a closed braid.
All of the genus one concordance operations are described by attaching two oriented
bands to the knot $K$, or rather to $K\times I$ = a neighborhood of $K$
in an orientable spanning surface for $K$; see Figure \ref{fig:genus1}. If the attaching 
regions of the bands represent linked pairs of points around $K$, the
result of the two band attachments is a twice punctured torus that
represents a genus one concordance between $K$ and another knot.

\begin{figure}[h]
\begin{center}
\includegraphics[width=3.5in]{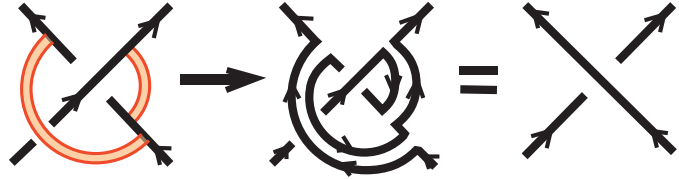}
\caption{Bands for a crossing change}\label{fig:cross}
\end{center}
\end{figure}

\begin{figure}[h]
\begin{center}
\includegraphics[width=4.5in]{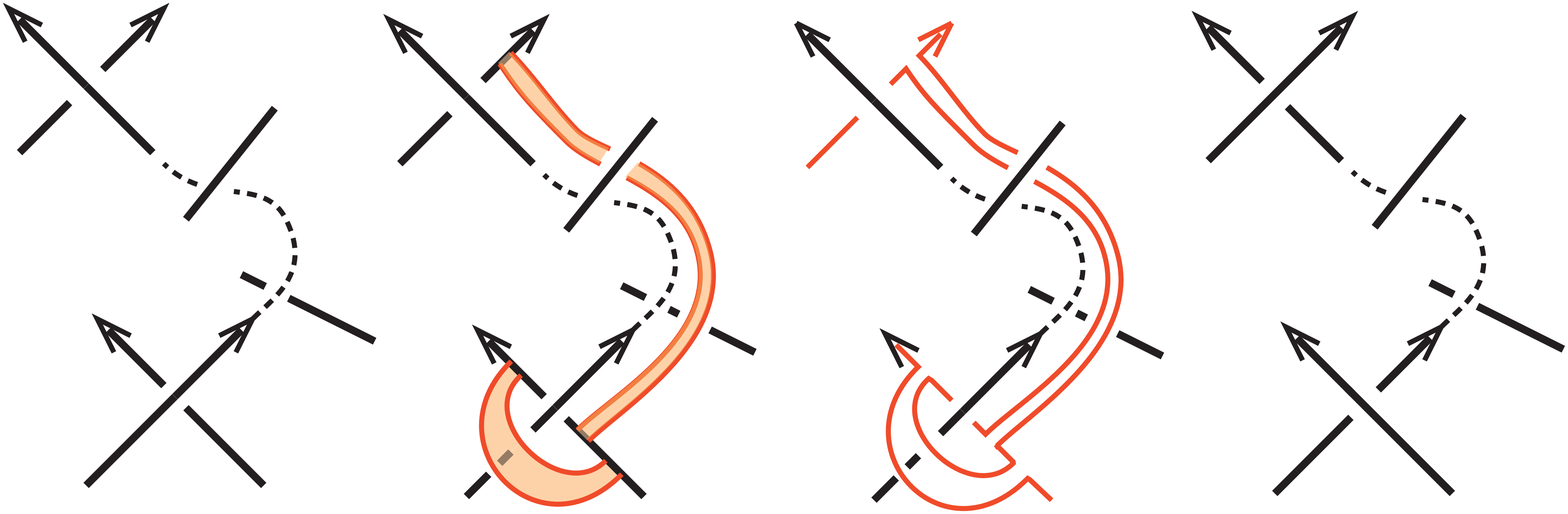}
\caption{Bands for switching positive and negative crossings}\label{fig:switch}
\end{center}
\end{figure}

\begin{figure}[h]
\begin{center}
\includegraphics[width=2in]{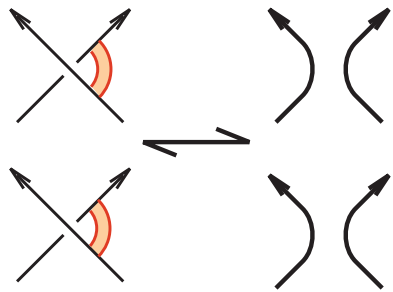}
\caption{Bands for resolving two crossings}\label{fig:resolve}
\end{center}
\end{figure}

The first of these genus one concordance operations is {\it crossing change}; that is,
if $K_1$ can be obtained from $K_0$ by a crossing change, then 
there is a genus one concordance between the two knots. This
is because a crossing change can be achieved by attaching two oriented
bands to $K_0$ (see Figure \ref{fig:cross}). 
The second operation, which we refer to as {\it switching two crossings}, consists 
of simultaneously switching a positive and a negative crossing. This
gives rise to genus one concordance since this again can be achieved
by attaching a pair of oriented bands (Figure \ref{fig:switch}). 
The third is {\it resolving two 
crossings} in a diagram for $K_0$; that is, removing the crossings so that 
an orientation on $K_0$ is preserved. This again can be achieved by attaching 
two bands, and hence, so long as the resulting link is again a knot $K_1$, 
there is a genus one concordance between the two knots (Figure \ref{fig:resolve}). 
The reverse of the resolving process, {\it de-resolving two crossings}
(in which two crossings are inserted), 
can similarly be achieved by attaching two bands by running the process in 
Figure \ref{fig:resolve} in reverse; this is the fourth of the
genus one concordance operations we employed in our computational search.

\section{First step: A braid representation data set and braid operations}\label{sec:braid}


For a knot $K_0$ given as the closure of a braid, the
four genus one concordance operations in Section~\ref{sec:genus1} can be 
carried out directly from the braid word representation of the knot. 
Denoting the standard generators for a
braid on $n$ strands by $\sigma_1,\ldots, \sigma_{n-1}$, 
where $\sigma_i$ represents a crossing of the $i$-th
strand over the $(i+1)$-st strand (reading upward), we use a common
shorthand that $i$ stands for $\sigma_i$ and $-i$ stands for $\sigma_i^{-1}$.
With this notation, a crossing change is effected by changing $i$ to 
$-i$, and simultaneously switching a positive and a negative crossing
is effected by changing a positive entry to its negative and a negative
entry to its positive in the braid word.
Resolving two crossings is effected by deleting
two generators (or their inverses) from the braid word,
and the reverse, or de-resolution, is effected 
by inserting two generators or their inverses into the braid word.
Note that in a knot diagram for $K_0$ that is the closure of a braid,
all of the strands in the braid can be oriented in the same direction (e.g., upward
in the braid),
and inserting or deleting any two braid generators preserves the orientation. 
All four of these braid word operations therefore performs one of
the genus one concordance operations in Section~\ref{sec:genus1} and changes the 4-genus 
of the underlying knot by at most one.

The first step in our process for determining 4-genus, and for 
determining genus one concordances with knots of known 4-genus
via the four crossing
operations, 
is to generate a database of representations of
the 21 knots $K_0$ listed in Theorem~\ref{theorem:slice} as closures of
braids.  
In fact, braid representatives of these 21 knots were found using a larger database
of randomly-generated braids which we had built for 
an ongoing project with a goal of determining the unknotting number of some of the
(as of this writing) 668 knots which the KnotInfo site reports have
unknown unknotting number; this database was also
used in our disproof of the Bernhard-Jablan Conjecture~\cite{britherm21}.

The construction of this larger database starts with a procedure
for generating random braids on 4 to 12 strands, using random braid
words of length between 20 and 60.
A python script is used to check that the
closure of the braid is a knot, and if so, the braid is added to the database
along with the volume of the knot computed using SnapPy~\cite{snappy}.
For each of the 21 knots $K_0$ in Theorem~\ref{theorem:slice}, 
and for each of the closed braid diagrams for knots $K_1$ in this larger database 
 whose volume matched the volume of $K_0$ to
within several decimal places,
we used SnapPy (in particular the command ``M0.is\_isometric\_to(M1)''
for the knot exteriors M0 and M1 of $K_0$ and $K_1$) to
check directly whether the knots were equal;
if so, the
braid closure representation was added to the database for that knot $K_0$.
This gave us our initial collection of braid representatives for the 21 knots.


\section{Pulling 4-genus down}\label{sec:down}


Using the braids in the database for the 21 knots $K_0$
discussed in Section~\ref{sec:braid}, 
we carried out all possible instances of
the three braid word operations corresponding to the
first three genus one concordance operations (crossing change,
switching two crossings, and crossing resolution), namely changing the sign of
a single letter, simultaneously changing a positive entry to its negative and a negative
entry to its positive, and deleting
two letters (or their inverses), in the braid word.
For each braid we considered from the database,
we also carried out 800 random de-resolutions, by randomly generated insertions
of two letters.

For each resulting braid, we then 
check whether the knot $K_1$ obtained by taking the closure of the braid
satisfies $g_4(K_1) = 0$ by using SnapPy to identify the knot $K_1$
and then checking whether that knot is in the list
from KnotInfo of knots (up to 12 crossings) that are smoothly slice.
We note that this part of the project was done before the announcement
on the arXiv (in April 2021) of the results of Karageorghis and Swenton~\cite{karswe21}
determining the 4-genera of the 18 knots in Theorem~\ref{theorem:slice}(b),
and at that time it was known, from the KnotInfo
database, that all 21 of the knots $K_0$ have smooth 4-genus either 1 or 2.
Thus a genus one concordance with a smoothly slice knot `pulled 
the 4-genus down', to determine that $g_4(K_0)=1$.
For the 18 knots in Theorem~\ref{theorem:slice}(b), our efforts 
succeeded. In particular, for 17 of these 18 knots, we have
a genus one concordance between the knot $K_0$, viewed as a braid closure, and a 
smoothly slice knot using a resolution of two crossings,
and for the 18th knot, we have a genus one concordance to a smoothly slice knot
constructed by switching two crossings.

In full detail, for each of the 18 knots we list a braid representative for it, the
slice knot (and the transformed braid representative) that is genus one concordant to it,
and the type of genus one concordance which yields the 4-genus determination. 
For every knot, there are,
in the end, many braids and genus one concordances to choose from; 
we have selected from among the
shortest braid representatives found.

\medskip

\begin{tabular}{ll}

$K11n80$ & $[5,3,-4,2,-4,6,-1,2,-3,-1,-1,4,2,1,6,5,6,-5,2,3,-6,5,4,-2]$\\
$6_1$ & $[5,\phantom{3}-4,2,-4,6,-1,2,-3,-1,-1,4,\phantom{2,,}1,6,5,6,-5,2,3,-6,5,4,-2]$\\
&(resolve two crossings)\\
\\

$K12a187$ & 
$[2,-3,4,6,5,5,-2,-6,-1,4,-3,-6,4,2,2,5,3,4,-3,4,-3,-5]$\\
$6_1$ & 
$[2,-3,4,6,5,5,-2,-6,-1,\phantom{4,}-3,-6,4,2,2,5,3,\phantom{4}-3,4,-3,-5]$\\
&(resolve two crossings)\\
\\

$K12a230$ & 
$[-4,3,-2,1,2,-3,4,4,-5,6,-5,3,6,4,2,3,2,-1,4,-5,-3,-3]$\\
$8_{20}$ &
$[-4,\phantom{3}-2,1,2,-3,4,4,-5,\phantom{6,}-5,3,6,4,2,3,2,-1,4,-5,-3,-3]$\\
&(resolve two crossings)\\
\\

$K12a317$ &
$[4,5,5,-3,-2,-3,4,4,3,5,-4,3,-2,3,-5,1,2,2,-3]$\\
$6_1$ &
$[4,5,\phantom{5,}-3,-2,-3,4,4,3,5,-4,3,-2,\phantom{3}-5,1,2,2,-3]$\\
&(resolve two crossings)\\
\\

$K12a450$ &
$[1,-4,-3,1,2,-6,4,-3,4,-6,2,3,5,-4,5,-1,-3,2,6,4]$\\
$6_1$ &
$[1,-4,-3,1,2,-6,4,-3,\phantom{4,}-6,2,3,5,-4,5,-1,-3,\phantom{2,}6,4]$\\
&(resolve two crossings)\\
\\

$K12a570$ &
$[-4,6,6,-4,3,5,4,2,-1,-1,-2,5,-6,-4,5,3,-2,3,5,2,1,2,4,-5]$\\
$6_1$ &
$[-4,\phantom{6,}6,-4,3,5,4,\phantom{2,}-1,-1,-2,5,-6,-4,5,3,-2,3,5,2,1,2,4,-5]$\\
&(resolve two crossings)\\
\\

$K12a624$ &
$[3,-2,6,-5,-6,-4,-3,5,5,6,-5,-5,-6,1,5,2,4,6,3,-5,4,5,2,2]$\\
$6_1$ &
$[3,-2,6,-5,-6,-4,-3,5,5,\phantom{6,}-5,-5,-6,1,5,\phantom{2,}4,6,3,-5,4,5,2,2]$\\
&(resolve two crossings)\\
\\

\end{tabular}

\begin{tabular}{ll}

$K12a636$ &
$[-3,-6,2,-1,-2,-5,4,5,-3,3,6,6,5,4,-2,3,-2,6,1,-5,3,4]$\\
$6_1$ &
$[-3,-6,2,-1,-2,-5,4,5,-3,\phantom{3,}6,6,5,4,-2,3,-2,6,1,-5,3\phantom{,4}]$\\
&(resolve two crossings)\\
\\

$K12a905$ & 
$[1,2,{\underline{2}},-3,2,-3,-3,-4,-2,3,4,4,-3,-2,1,3,2,{\underline{-2}},3,-2]$\\
$10_{87}$ &
$[1,2,{\underline{-2}},-3,2,-3,-3,-4,-2,3,4,4,-3,-2,1,3,2,{\underline{2}},3,-2]$\\
&(switch two crossings, marked with {\underline{underlining}})\\
\\

$K12a1189$ & 
$[2,2,-4,3,3,-2,4,1,2,1,-3,-4,2,-1,2,-3,-2,1,-4,3]$\\
$10_{137}$ &
$[2,2,-4,3,3,-2,4,1,2,\phantom{1,}-3,-4,2,-1,2,-3,-2,\phantom{1,}-4,3]$\\
&(resolve two crossings)\\
\\

$K12a1208$ & 
$[-2,4,1,-3,4,-3,4,-3,2,-4,3,3,-2,1,1,4,-2,3,3,-1,4,-3]$\\
$10_{137}$ &
$[-2,4,1,-3,4,-3,\phantom{4,}-3,2,-4,3,3,-2,1,1,4,-2,3,3,-1,\phantom{4,}-3]$\\
&(resolve two crossings)\\
\\

$K12n52$ &
$[3,-1,6,4,-1,-5,4,6,-3,2,-5,2,1,3,5,-5,-5,4,5,3,-4,2]$\\
$10_{75}$ &
$[3,-1,6,4,-1,-5,\phantom{4,}6,-3,2,-5,2,1,3,\phantom{5,}-5,-5,4,5,3,-4,2]$\\
&(resolve two crossings)\\
\\

$K12n63$ &
$[-5,-5,3,-2,3,3,2,6,2,-1,-3,2,4,-1,5,6,6,3,2,-3,-1,-3,4,1]$\\
$6_1$ &
$[-5,-5,3,-2,3,3,2,6,2,-1,-3,2,\phantom{4,}-1,5,6,6,3,2,-3,-1,-3,4\phantom{,1}]$\\
&(resolve two crossings)\\
\\

$K12n225$ &
$[-2,5,1,-3,-6,-6,4,6,-6,1,-2,4,4,-5,-3,2,-6,2,-4,1,-2,5,4,4]$\\
$K11n74$ &
$[-2,5,1,-3,-6,-6,4,\phantom{6,}-6,1,-2,4,4,-5,-3,2,-6,2,-4,1,-2,\phantom{5,}4,4]$\\
&(resolve two crossings)\\
\\

$K12n555$ &
$[5,-3,-6,5,7,4,6,1,5,-3,4,-2,3,5,4,-6,-3,-6,-3,4,7]$\\
$6_1$ &
$[5,-3,-6,\phantom{5,}7,4,6,1,5,-3,4,-2,3,5,4,-6,-3,-6,-3,\phantom{4,}7]$\\
&(resolve two crossings)\\
\\

\end{tabular}

\begin{tabular}{ll}

$K12n558$ &
$[3,-4,-1,-1,3,-1,4,-2,3,1,-4,3,1,4,-3,2,3,2]$\\
$8_8$ &
$[3,-4,-1,-1,3,-1,4,-2,3,\phantom{1,}-4,3,1,4,-3,2,3\phantom{,2}]$\\
&(resolve two crossings)\\
\\

$K12n665$ &
$[1,6,4,5,-5,-1,4,-3,5,2,-3,-1,2,4,-5,3,3,-2,-5,3,2,1,-2,3,-2,3]$\\
$K12n256$ &
$[1,6,\phantom{4,5,}-5,-1,4,-3,5,2,-3,-1,2,4,-5,3,3,-2,-5,3,2,1,-2,3,-2,3]$\\
&(resolve two crossings)\\
\\

$K12n886$ &
$[4,2,-3,1,-5,4,-3,-2,3,4,3,5,2,2,4,4,-3,-3,-1,-1,-2,-5,-3,2,3]$\\
$6_1$ &
$[4,2,-3,1,-5,4,-3,-2,3,4,3,5,2,2,4,\phantom{4,}-3,-3,-1,-1,-2,-5,\phantom{-3,}2,3]$\\
&(resolve two crossings)\\
\\

\end{tabular}

\medskip

This completes the proof of Theorem~\ref{theorem:slice}(b).
As illustrated by the above list, the genus one concordance
operation of resolving two crossings has been by 
far the most successful strategy.


\section{Pulling 4-genus up}\label{sec:up}

For the remaining three knots $K_0$ with unknown smooth 4-genus, $K12a153$, 
$K12n239$, and $K12n512$, the procedure in Section~\ref{sec:down}
failed to find a slice knot that is genus one concordant with $K_0$ after running 
on 1800, 20000, and 3500 closed braid diagrams
of the knot, respectively (and after running for several months, 
not including the time for creating the database).
On the other hand, that procedure succeeded in finding the genus one concordances for
the 18 knots discussed in Section~\ref{sec:down} after running on 1 - 500 diagrams
for each knot, and on less than 10 for most of them (and within a few days).
This naturally led to the
conclusion that, since all three knots are known to have 4-genus either 1 or 2, they
most likely have smooth 4-genus equal to~2. 

In the same computer search described in Section~\ref{sec:down}, in which
the four braid word operations were performed on
braid representations in our database for the three knots $K12a153$, 
$K12n239$, and $K12n512$, for each resulting braid we also 
checked whether the knot $K_1$ obtained by taking the closure of the braid
satisfies $g_4(K_1) = 3$. 
As above this process uses SnapPy to identify the knot $K_1$
and then checks whether that knot is in the ``target list''
consisting of the set of knots from KnotInfo that have smooth 4-genus equal to~3.
This attempt to pull the 4-genus up also did not succeed.
However, in this case the situation is complicated by the fact that
all three of the knots $K12a153$, 
$K12n239$, and $K12n512$ have topological 4-genus equal to~1,
according to the KnotInfo site. In analogy with Equation~(\ref{eq:concordance}),
the topological 4-genera of knots that are genus one concordant differ by at most~1
(by an identical argument to the one in Section~\ref{sec:genus1} for smooth 
4-genus),
and so any knot sharing a 
genus one concordance with one of the three knots must have topological
4-genus at most 2. There are only 24 knots with 12 or fewer
crossings which have smooth 4-genus equal to 3 and also topological
4-genus at most 2 (again from the KnotInfo site), and thus
effectively the number of potential knots in the ``target list'' was
limited to only 24.

As a consequence, we
created a second procedure for pulling up the 4-genus for the remaining
three knots, with modifications both to increase the number of diagrams considered 
for each of the three knots and to increase the number of knots in the target list.

We describe the method of increasing the target list first.
In order to construct a larger list of
knots with smooth 4-genus (at least) 3 
we adopt another result from the work of Lewark and McCoy,
which is an application of Donaldson's 
Diagonalization Theorem \cite{don87}. 
We restate 
their result here.

\begin{theorem}\cite[Lemma~3]{lewmcc19}\label{theorem:lmt}
Let $K$ be a knot with a positive-definite $m \times m$ Goeritz matrix $G$.
If $\sigma(K)\leq 0$ and $2g_4(K) = -\sigma(K)$, then there is an $(m-\sigma(K))\times m$ 
integer matrix $M$, such that $G = M^TM$.
\end{theorem}

Lewark and McCoy use the contrapositive of this result to show
that twelve 11- and 12-crossing 
knots $K$ with previously unknown smooth 4-genus and with $\sigma(K)=-2$
satisfy $g_4(K)\geq 2$ (and hence $g_4(K) = 2$). 
We also apply the contrapositive of Theorem~\ref{theorem:lmt} to build our
larger target list of knots with 4-genus at least 3, as follows.


Using SnapPy's census of knots with crossing number 13 through 16,
contained in the
``AlternatingKnotExteriors'' and
``NonalternatingKnotExteriors'' databases, as our source data, 
we used SnapPy to find the knots $K$ with signature equal to~-4.
We note that some care needs to be taken in implementing this step, since one can find
different (contradictory) conventions on the sign of the signature of a knot,
and negative signature plays a very important role in Lewark and McCoy's
result. In particular, SnapPy~\cite{snappy}, prior to version 3.0, implemented the `opposite'
sign convention for the signature; this was changed, however, in April 2021 
with version 3.0.
Among the knots with signature~-4, we then
used SnapPy in Sage~\cite{sage} to carry out matrix factorization
computations, via Sage's built-in implementation of GAP \cite{GAP4},
in order to find the knots $K$ for which the Goeritz matrix
does not admit the decomposition in Theorem~\ref{theorem:lmt},
and hence for which $g_4(K) \ge 3$.

Our implementation of this search found no knots $K$ with 14 or fewer crossings,
with $\sigma(K)=-4$, for which Theorem~\ref{theorem:lmt} could conclude that $g_4(K)\geq 3$.
We did, however, find three 15-crossing alternating knots and 27 
16-crossing alternating knots
for which we could conclude that $g_4(K)\geq 3$. 

We also searched for target knots among the 1,769,978 hyperbolic
17-crossing alternating knots, using
the DT codes contained in the tables of knots provided with the software package
Regina~\cite{bur20}.
This yielded~332 alternating 17-crossing knots with signature equal to~$-4$ which we could 
show have smooth 4-genus
at least 3. Using the 6,283,385 hyperbolic non-alternating 17-crossing knots gave a further 49
examples. This list of $24 + 3 + 27 + 332 + 49 = 435$ knots of crossing number up to 17
and 4-genus at least 3 is the target list of our second procedure.

In our second procedure we employed two methods to build a larger data set of
diagrams for $K_0 \in \{K12a153, K12n239, K12n512\}$,
and of diagrams for knots $K_1$ that are genus one concordant to $K_0$.

The first method starts with a minimal crossing diagram for $K_0$,
uses the ``backtrack'' function in SnapPy to randomaly generate another diagram
for $K_0$ with more crossings, and then applies the ``braid\_word'' function in SnapPy
to output a representation of $K_0$ as a closure of a braid.
Then as before, this method generates diagrams for knots $K_1$ genus one concordant to $K_0$ by
carrying out all possible instances of
the braid word operations corresponding to crossing change,
switching two crossings, and crossing resolution, along with 800 randomly
generated insertions of two letters corresponding to de-resolutions.

The second method also starts with a minimal crossing diagram for $K_0$ and
uses ``backtrack''to randomaly generate a new diagram
for $K_0$ with more crossings, and then uses SnapPy's ``PD\_code'' to construct
a planar diagram (PD) code for that diagram.
We have written a separate algorithm in Python that carries
out operations on PD codes corresponding to crossing change,
switching two crossings, and resolving two crossings (although not de-resolution), 
and this is carried out for all possible operations on the PD code of the diagram for $K_0$
to create diagrams for knots $K_1$ genus one concordant to $K_0$.

The last step of the second procedure checks whether a knot $K_1$ in our data set,
that is genus one concordant to $K_0$, 
is equal to one of the 435 knots in our target list of knots with 4-genus at least 3.
Unlike our first procedure in Section~\ref{sec:down}, the diagrams involved
here have more crossings. In this case, the procedure uses SnapPy to compute 
the volume of the complement of $K_1$ and then to check whether that complement
is isometric to the complements of any of the 435 target knots whose volume
matches that of $S^3 \setminus K_1$ to several decimal places.

\medskip

\begin{figure}[h] \label{fig:K153}
\begin{center}
\includegraphics[width=3.5in]{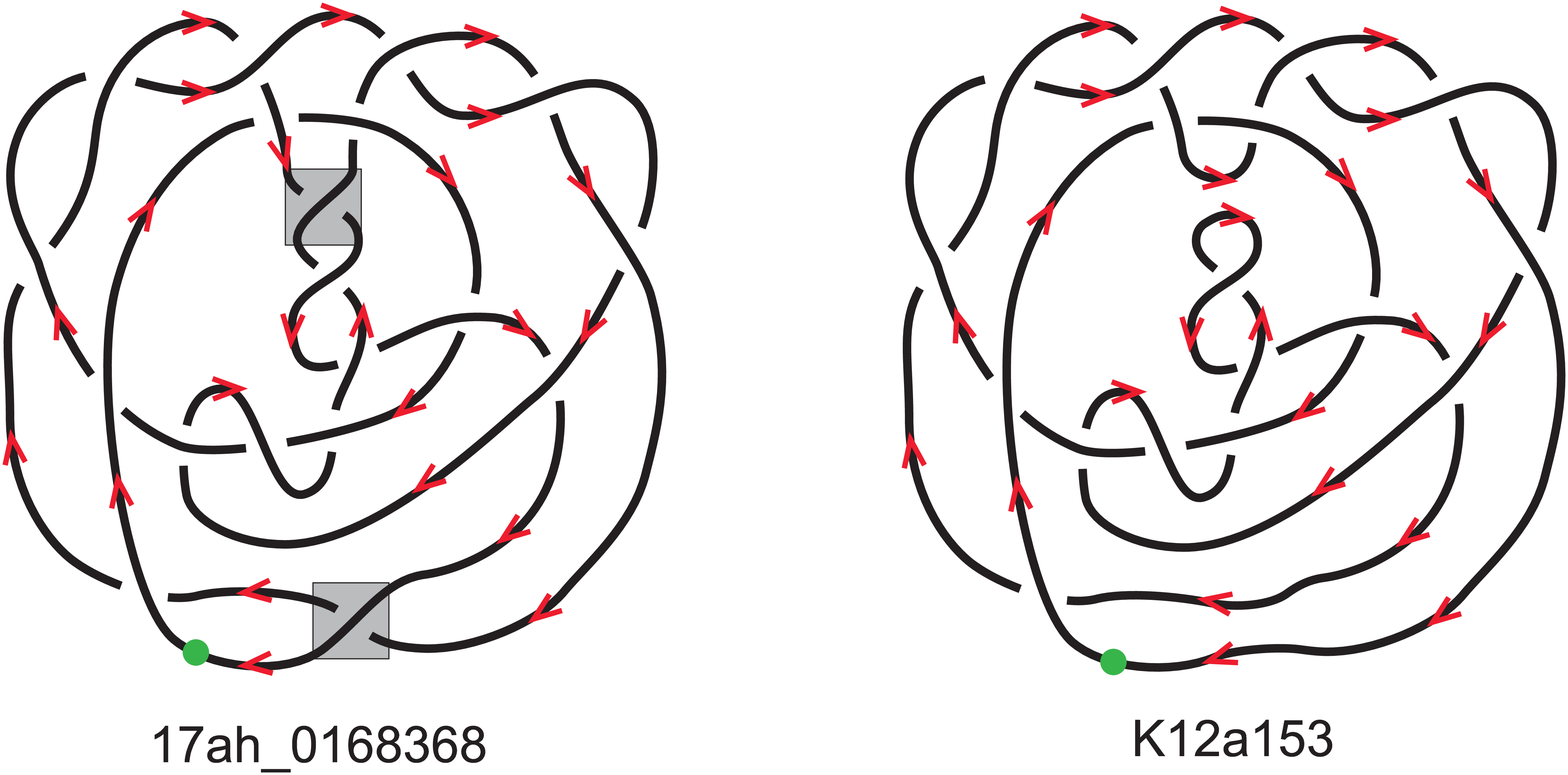}
\caption{A (de-)resolution for K12a153}\label{fig:K12a153}
\end{center}
\end{figure}

\begin{figure}[h] \label{fig:K239}
\begin{center}
\includegraphics[width=3.5in]{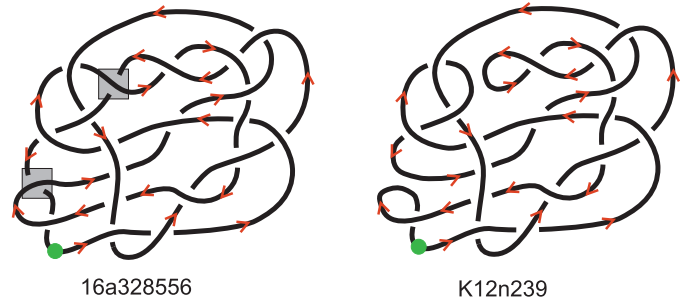}
\caption{A (de-)resolution for K12n239}\label{fig:K12n239}
\end{center}
\end{figure}

For the knots $K12a153$ and $K12n239$, both the planar diagram and the braid closure
approaches succeeded. The planar diagram approach found smaller diagrams in general, and 
in particular, after further simplification (via Reidemeister moves), 
Figures~\ref{fig:K12a153} and~\ref{fig:K12n239} illustrate 
17-crossing and 18-crossing diagrams for
the knots $K12a153$ and $K12n239$ together with de-resolutions of two crossings 
yielding the knots $17ah\_0168368$ (from Regina's database)
and $16a328556$ (from SnapPy's database), respectively. 
A computation shows that the knots $17ah\_0168368$
and $16a328556$ meet the conditions of the contrapositive of
Theorem~\ref{theorem:lmt}; their Goeritz matrices can be found in an appendix at the end of the
paper. These knots therefore have 4-genera at least 3, and hence
(from Equation~\ref{eq:concordance} and the fact that the 4-genus is either 
1 or 2 from KnotInfo) the knots $K12a153$ and $K12n239$
both have 4-genus equal to 2.

DT codes for these knot diagrams (starting from the indicated points) are:

\smallskip

\noindent $K12a153$ :\hfill

\noindent $[26,-30,-24,-18,-2,-28,-32,-34,-8,-6,-22,-20,-16,-10,-12,-4,-14]$\hfill


\noindent $17ah\_0168368$ :\hfill

\noindent $[28,-32,-36,-18,-2,-30,-24,-26,-8,-6,-34,-4,-14,-38,-10,-12,-22,-20,-16]$\hfill

\smallskip

\noindent $K12n239$ :\hfill

\noindent $[20,24,12,10,16,-18,30,8,6,34,2,32,14,28,26,4,22,36]$\hfill



\noindent $16a328556$ :\hfill

\noindent $[28,32,20,38,36,34,4,30,40,-26,12,8,6,16,2,14,22,10,24,-18]$\hfill


\medskip

On the other hand, the knot $K12n512$ proved more difficult to capture by these methods. In particular, 
after several months of searching, none of the 435 `target' knots 
with 4-genus at least 3 that we had identified
using Theorem \ref{theorem:lmt} managed to turn up
in our data of knots genus one concordant to $K12n512$. 
We then expanded our target list
of knots $K$ with $g_4(K)\geq 3$, by searching through
the 8,400,285 alternating hyperbolic 18-crossing knots from the Regina database
and again using 
Theorem~\ref{theorem:lmt}, 
yielding 3544 further examples for our target list. 
(Out of 1,958,354 knots with signature $\pm 4$, this was a 
$0.18\%$ `success' rate). Running searches on 
the knots genus one concordant to the knot $K12n512$, built
from roughly 200,000
diagrams for $K12n512$,
succeeded using the braid closure approach, using the operation of
de-resolution of two crossings.
In particular, our procedure yielded the 11-strand braid

\smallskip

$[4, 6, 3, -5, 2, -4, 1, 3, 7, 4, -6, 7, -6, -8, -6, 9, 5, -10, -6, 9, -5, 7, -4, -3, 8, 2,$

$ -7, -6, -5, 4, 3, 5, 7, 4, -6, 5, -7, -8, -7, -9, -6, 10, -7, -8, 7, 6, -5, -4, -3,$

$ -2, -1, -2, -3, -4, 5, 6, -7, 8, -9, 8, -7, -6, -5, -7]$

\smallskip

\nidt for the knot $K12n512$, and the braid

\smallskip

$[4, 6, 3, -5, 2, -4, 1, 3, 7, 4, -6, 7, -6, -8, -6, 9, 5, -10, -6, 9, -5, 7, -4, -3, 8, 2,$

$ -7, -6, -5, 4, 3, 5, 7, 4, -6, 5, -7, -8, -7, -9, -6, 10, -7, -8, 7, 6, -5, -4, -3,$

$ -2, -1, -2, -3, -4, 5, 6, -7, 8, -9, 8, -7, -6, -5, -7, -8, -1]$

\smallskip

\nidt for the knot $K_1 = 18ah\_2335674$ (from the Regina database),
which has $g_4(K_1)\geq 3$. 
The Goeritz matrix for the latter knot
can also be found in the appendix. The knot 
$18ah\_2335674$ is obtained from $K12n512$ by inserting
the last two entries of the braid word for $18ah\_2335674$, 
representing the de-resolution of the two crossings.
Consequently, $K12n512$ has smooth 4-genus at least 2, and so (with the information 
that the genus is at most 2 from KnotInfo) we have
$g_4(K12n512)=2$.

With this result, the proof of Theorem~\ref{theorem:slice} is complete; 
the knots $K12a153$, $K12n239$ and $K12n512$
all have smooth 4-genus equal to 2. This finishes the computation of smooth 4-genus for 
all prime knots though 12 crossings.


\section{Where to go from here?}\label{sec:future}


Although the smooth 4-genus of every prime knot with 12 or fewer crossings
is now known, this still does not allow us to immediately 
compute the 4-genus of the non-prime knots in the same crossing range. The 3-genus is well-known 
to be additive under connected sum \cite{schu}, but it is equally well-known that
the 4-genus is not; for example, the connected sum $K\sharp(-K)$ of a knot with its
mirror image is (ribbon, and hence) a slice knot, so has smooth 4-genus zero. 
The computation of $g_4(K)$ for connected sums of knots is a challenging problem, and 
is closely related to the idea of
cobordism distance~\cite{cohar18},\cite{felpar21}.
There are several interesting results on this problem, especially for connected sums of 
torus knots \cite{liv18},\cite{fel16}, but it is far from complete.

Turning to {\it topological} 4-genus, there still remain, according to KnotInfo, 7 prime
knots through 12 crossings whose topological 4-genus is unknown, namely
$K12a244$, $K12a810$, $K12a905$, $K12a1142$, $K12n549$, $K12n555$, and $K12n642$. 
All have $g_4^{top}(K)$ equal to $1$ or $2$. It is possible that techniques like the
ones carried out here could shed light on these remaining cases, as well; this
is something that we plan to explore in the future.

The notion of a {\it concordance} (of genus zero) - that is,
a (smooth) annulus in $S^3\times I$  between knots $K_0$ and $K_1$ - 
is at the heart of still another knot invariant which 
essentially sits between 
$g_4(K)$ and $g_3(K)$: the concordance genus $g_c(K)$ is the minimum 3-genus of a
knot concordant to $K$~\cite{liv04}. Its value is known for every prime knot through 10
crossings, but its value is unknown (at present) for 208 prime knots through 
12 crossings. It would be interesting to see if techniques similar to those 
employed here could make progress on shrinking this list.


\section{Appendix: matrix computations}\label{sec:append}

Here we provide Goeritz matrices for the three knots listed in Section~\ref{sec:up}
found to have
$g_4(K)\geq 3$ and used to pull the smooth 4-genera of the knots
$K12a153$, $K12n239$ and $K12n512$ up to 2. 

\bigskip

$17ah\_0168368$ : DT code = $[18,20,30,28,24,4,8,26,10,2,32,34,16,14,6,12,22]$

\smallskip

Goeritz matrix  = $\left(
\begin{array}{@{}*{11}{c}@{}}
5&-3&-1&-1&0&0&0&0&0\\
-3&6&0&0&-3&0&0&0&0\\
-1&0&3&0&-1&-1&0&0&0\\
-1&0&0&3&-1&0&-1&0&0\\
0&-3&-1&-1&6&0&0&-1&0\\
0&0&-1&0&0&2&0&0&-1\\
0&0&0&-1&0&0&2&0&0\\
0&0&0&0&-1&0&0&3&-1\\
0&0&0&0&0&-1&0&-1&2
\end{array}
\right)$

\bigskip

$K16a328556$ : DT code = $[6, 14, 26, 16, 20, 22, 4, 24, 32, 30, 28, 12, 2, 10, 8, 18]$

\smallskip

Goeritz matrix = $\left(
\begin{array}{@{}*{11}{c}@{}}
3&-1&-1&-1&0&0&0&0&0\\
-1&3&0&0&-1&-1&0&0&0\\
-1&0&3&0&0&-1&-1&0&0\\
-1&0&0&5&0&-2&-1&-1&0\\
0&-1&0&0&2&0&0&0&-1\\
0&-1&-1&-2&0&6&0&-1&0\\
0&0&-1&-1&0&0&2&0&0\\
0&0&0&-1&0&-1&0&3&0\\
0&0&0&0&-1&0&0&0&2
\end{array}
\right)$

\bigskip

$18ah\_2335674$ : DT code = $[32,14,30,28,24,22,36,4,6,26,12,10,34,18,16,2,20,8]$

\smallskip

Goeritz matrix = 
$\left(
\begin{array}{@{}*{11}{c}@{}}
2&-1&-1&0&0&0&0&0&0&0&0\\
-1&4&0&-1&-1&-1&0&0&0&0&0\\
-1&0&4&-2&0&0&-1&0&0&0&0\\
0&-1&-2&6&0&0&-1&-1&-1&0&0\\
0&-1&0&0&2&0&0&0&0&-1&0\\
0&-1&0&0&0&3&0&0&-1&0&-1\\
0&0&-1&-1&0&0&3&-1&0&0&0\\
0&0&0&-1&0&0&-1&3&0&0&-1\\
0&0&0&-1&0&-1&0&0&3&0&0\\
0&0&0&0&-1&0&0&0&0&2&0\\
0&0&0&0&0&-1&0&-1&0&0&2
\end{array}
\right)$

\medskip

In addition, files containing the names and DT codes of all of the knots
that we identified using Theorem \ref{theorem:lmt} to have 4-genus at least $3$ and 
signature $-4$ can be found on the authors' website, at

\begin{center}
https://www.math.unl.edu/$\sim${mbrittenham2}/{knot}\_data/{4genus}/
\end{center}

\noindent along with the python code used to find the knots.







\begin{thebibliography}{99}


\bibitem{bur20}
B. A. Burton,
{\it The next 350 million knots}, 
36th International Symposium on Computational Geometry (SoCG 2020) (S. Cabello, D.Z. Chen, eds.), Leibniz Int. Proc. Inform.
{\bf 164} 
Schloss Dagstuhl–Leibniz-Zentrum fuer Informatik
(2020) 
25:1--25:17.

\bibitem{britherm21}
M. Brittenham and S. Hermiller,
{\it A counterexample to the Bernhard–Jablan Unknotting Conjecture},
Exp. Math. {\bf 30} (2021), no. 4, 547--556.

\bibitem{cohar18}
T. Cochran and S. Harvey,
{\it The geometry of the knot concordance space},
Algebr. Geom. Topol. {\bf 18} (2018) 2509--2540.

\bibitem{don87}
S. K. Donaldson, 
{\it The orientation of Yang-Mills moduli spaces and 4-manifold topology}, 
J. Diff. Geom. {\bf 26} (1987) 397--428.

\bibitem{fel16}
P. Feller, {\it Optimal Cobordisms between Torus Knots},
Comm. Anal. Geom. {\bf 24} (2016) 993--1025.

\bibitem{felpar21}
P. Feller and J. Park,
{\it Genus One Cobordisms Between Torus Knots},
Inter. Math. Res. Notices {\bf 2021} (2021) 521--548.

\bibitem{GAP4}
The GAP~Group,
\newblock {\em {GAP -- Groups, Algorithms, and Programming, Version 4.7.8}},
2015.

\bibitem{juh10}
A. Juh\'asz, 
{\it The sutured Floer homology polytope}, 
Geom. Topol. {\bf 14} (2010) 1303--1354.

\bibitem{karswe21} L. P. Karageorghis and F. Swenton, 
{\it Determining the doubly slice genera of prime knots with up to 12 crossings},
J. Knot Theory Ramifications {\bf 30} (2021), no. 8, Paper No. 2150057, 17 pp.

\bibitem{km97}
P. Kronheimer and T. Mrowka, 
{\it Scalar curvature and the Thurston norm},
Math. Res. Lett. {\bf 4} (1997) 931--937.



\bibitem{lewmcc19}
L. Lewark and D. McCoy,
{\it On calculating the slice genera of 11- and 12-crossing knots},
Experimental Math. {\bf 28} (2019) 81--94.

\bibitem{liv04}
C. Livingston, {\it  The concordance genus of knots},
Algebr. Geom. Topol. {\bf 4} (2004) 1--22.

\bibitem{knotinfo}
C. Livingston and A. Moore,
{\it KnotInfo: Table of Knot Invariants},
{\tt knotinfo.math.indiana.edu} (11/28/2021).

\bibitem{liv18}
C. Livingston and C. Van Cott, {\it The four-genus of connected sums of torus knots},
Math. Proc. Camb. Phil. Soc. {\bf 164} (2018) 531--550.

\bibitem{mur65}
K. Murasugi, {\it On a certain numerical invariant of link types},
Trans. Amer. Math. Soc. {\bf 117} (1965) 387--422.

\bibitem{oz04}
P. Ozsv\'ath and Z. Szab\'o,
{\it Holomorphic disks and knot invariants},
 Adv. Math. {\bf 186} (2004) 58--116.

\bibitem{oz04b}
P. Ozsv\'ath and Z. Szab\'o,
{\it Holomorphic disks and genus bounds}, Geom. Topol. {\bf 8} (2004) 311--334.

\bibitem{oz11}
 P. Ozsv\'ath and Z. Szab\'o, 
{\it Knot Floer homology and rational surgeries}, 
Alg. Geom. Top. {\bf 11} (2011) 1--68.

\bibitem{pic20}
L. Piccirillo,
{\it The Conway knot is not slice},
Annals of Mathematics
{\bf 191} (2020) 581--591.

\bibitem{ras03}
J. Rasmussen,
{\it Floer homology and knot complements}, PhD thesis, Harvard University
(2003).

\bibitem{ras10}
J. Rasmussen, {\it Khovanov homology and the slice genus},
Invent. Math. {\bf 182} (2010) 419--447.

\bibitem{sage}
The Sage Developers,
{\it {S}ageMath, the {S}age {M}athematics {S}oftware {S}ystem ({V}ersion 9.2)},
{\tt https://www.sagemath.org}, 2020.

\bibitem{schu}
H. Schubert, {\it Knoten und Vollringe}, Acta Mathematica 
{\bf 90} (1953) 131--286.

\bibitem{snappy}
M. Culler, N. Dunfield, M. Goerner and J. Weeks,
{\it Snap{P}y, a computer program for studying the geometry and topology of $3$-manifolds},
Available at {\tt http://snappy.computop.org} (11/28/2021).

\end{thebibliography}
\end{document}